\documentclass[12pt]{article}
\usepackage{amsmath,amsthm,amssymb,latexsym, amsfonts}
\usepackage[dvipsnames]{xcolor}
\usepackage[dvips]{graphicx}

\theoremstyle{plain}
\newtheorem{teor}{Theorem}[section]
\newtheorem{defin}[teor]{Definition}
\newtheorem{propo}[teor]{Proposition} 
\newtheorem{obs}[teor]{Remark} 
\newtheorem{theo}[teor]{Theorem}
 
\newtheorem{example}[teor]{Example} 
\newtheorem{coro}[teor]{Corollary} 

\def\proof{{\noindent \bf Proof:} \hspace{0.1 cm}}
\newcommand{\cqd}{\hfill \rule{2mm}{2mm}\vspace{.3cm}}

\def \Z {\mathbb Z}

\def \F {\mathbb F}

\date{}
\title{Fermat's Little Theorem and Euler's Theorem in a class of rings}

\begin{document}

\author{ Fernanda D. de Melo Hern\'andez \thanks{Corresponding author}
\thanks{F. D. de Melo Hernandez (fdmelo@uem.br), Departamento de Matem\'atica, Universidade Estadual de Maring\'a, Av. Colombo 5790, 87020-900, Maring\'a, PR, Brazil.} \and
C\'esar A. Hern\'andez Melo \thanks{C\'esar A. Hern\'andez Melo (cahmelo@uem.br), Departamento de Matem\'atica, Universidade Estadual de Maring\'a, Av. Colombo 5790, 87020-900, Maring\'a, PR, Brazil.} 
\and Horacio Tapia-Recillas \thanks{Horacio Tapia-Recillas (htr@xanum.uam.mx), Departamento de Matem\' aticas, Universidad Aut\' onoma Metropolitana-Iztapalapa, Av. San Rafael Atlixco 186, 09340, CDMX, M\' exico.}}
\maketitle


\begin{abstract}

Considering $\mathbb{Z}_n$ the ring of integers modulo $n$, the classical Fermat-Euler theorem establishes the existence of a specific natural number $\varphi(n)$ satisfying  the following property: 
\begin{equation}\label{one}
x^{\varphi(n)}=1
\end{equation}
for all $x$ belonging to the group of units of $\mathbb{Z}_n$. In this manuscript, this result is extended to a class of rings that satisfies some mild conditions. \\

\end{abstract}



\section{Introduction}

Fermat's Little Theorem (FLT) establishes that for a prime number $p$, $x^{p-1} \equiv 1$ mod $p$ for any integer $x$. For a positive integer $m$, Euler's Theorem (ET) establishes that $x^{\varphi(m)} \equiv 1$ mod $m$ for any integer $x$   relatively prime to $m$ where $\varphi$ is the Euler totient function given the number of integers relatively prime to $m$. If $m=p$, a prime, Fermat's Little Theorem is obtained from Euler's Theorem. These results have been studied for a long time, and nowadays have interesting applications in cryptography with the Diffie-Hellman Key Exchange, the ElGamal cryspotysterm, the public key cryptosystem RSA and related topics such as primality testing. The generalization of the Fermat and Euler Theorems have gone in several directions as it is shown in the references \cite{C}, \cite{D}, \cite{G}, \cite{HaSm}, \cite{IsPo}, \cite{M}, \cite{N-L}, \cite{Z}.\\

For instance in \cite{C} the Theorems of Wilson, Fermat and Euler are generalized by using Polya's and de Bruijn's theorems of enumeration.
In \cite{G} by introducing the concept of a Fermat sequences $\{u_n\}_{n=1}^{\infty}$ of rational numbers, i.e., a sequence having the property, $u_p \equiv u_1$ mod $p$, called a Fermat sequence, it si shown that there exits a class of Fermat sequences over the rational numbers.
 In \cite{IsPo} by using Group theory another proof os Fermat's Little Theorem is established.
 Certain rational integral functions on one indeterminate associated to the FLT is generalized to several indeterminates as it is shown in \cite{M}.\\

Given any prime $p$, FLT can be stated as the congruence $x^p - x \equiv 0$ (mod $p$) is satisfied by every integer $x$. Moreover, if $f(x)$ is a polynomial with integral coefficients such that $f(x)\equiv 0$ (mod $p$) is satisfied by every integer $x$, then $f(x)$ is a multiple of $x^p - x$. In \cite{N-L} the authors raise and settle the following question: for any positive integer $m$ characterize the class of polynomials $f(x)$ having the property that $f(x) \equiv 0$ (mod $m$) is satisfied by every integer. It turns out that in general, unlike the Fermat case where $m$ is a prime, these polynomials are not all multiples of a single polynomial. \\

For an integer $m$ an element $e \in \mathbb Z_m$ is called idempotent if $e^2 \equiv e $ mod (m). If $E_m$ denotes the set of idempotents of $\mathbb Z_m$, in \cite{V} among other results the following is proved: for all $m \in \mathbb N$, $a \in \mathbb Z_m$; $a^{\varphi(m)} \in E_m$. \\

In \cite{Z} several results related to the trace of matrices with entries in $\mathbb Z$ are proved, including the following (Arnold's conjecture): for any such a matrix $A$, for any prime $p$ and for any integer $a$: $Tr(A^n) \equiv Tr(A^{n-\varphi(n)})$ where $n=p^a$ and $\varphi(n)$ is the Euler function. \\

If $\mathbb Z_m$ is the ring of integers modulo $m$ and $\mathbb Z_m^{\ast}$ denotes the group of units of this ring, whose order is $\varphi(m)$, Euler's Theorem states that for any unit of $\mathbb Z_m$, 
\[
x^{|\mathbb Z_m^{\ast}|} \equiv 1 \; {\rm mod} \; m.
\]

In this manuscript a similar result to Euler Theorem is presented for a class of rings that includes the ring $\mathbb Z_m$ as a particular case. The paper is organized as follows:
in Section 2 basic results on the group of units of a ring, particularly the multiplicative order of a ring, are recalled. In Section 3 the main results of the manuscript are given (Theorems \ref{Fermat} and \ref{Euler}), and in Section 4 some application of these results are presented. Several examples are given illustrating the main results.


\section{Basic facts}\label{BF}

 Recall that an element $x$ of a ring $R$ with identity 1 is {\it invertible} if $y \in R$ exists such that $xy=1=yx$. The set of invertible elements of a ring $R$ will be denoted by $R^{\ast}$ and it is a group under the multiplication operation of the ring. Also, an element $x \in R$ has \textit{finite multiplicative order}, if there exists a natural number $m \geq 1$ such that  $x^m=1$. The smallest natural number different from zero $m$ with this property is called the \textit{multiplicative order} of the  element $x$ and it will be denoted by $o(x)$, i.e., 
\[
o(x)={\rm min}\{m \in \mathbb N: \; x^m=1  \}.
\]
It is clear that the identity of the ring $R$ is the unique element in the ring with multiplicative order equal to one. We assume that $0\in R$ has order 0. In the following lines basic results involving the concepts just introduced are presented.

\begin{itemize}
\item[1.]\label{f} Observe that the order of an element of a ring defines a function:
\[
O: R \longrightarrow \mathbb N, \; O(x)=o(x).  
\]

\item [2.]\label{u} If $x \in R\setminus\{0\}$ has finite multiplicative order then it is invertible.

\item[3.]\label{AAA} If $x$ has finite multiplicative order and $x^m=1$, then $o(x)$ divides $m$. In fact, applying the division algorithm in $\mathbb{Z}$, there exists $k,r\in\mathbb{Z}$ with $0\leq r< o(x)$ such that $m=o(x)k+r$, since
$$
1=x^m=x^{o(x)k+r}=(x^{o(x)})^kx^r=x^r,
$$
the definition of $o(x)$ implies that $r=0$, from where $o(x)|m$.

\item[4.] \label{Bb}  If $x\in R$ is an invertible element and it does not have finite multiplicative order, then the group of units $R^*$ of the ring $R$ is infinite. In fact, from the assumptions, it follows that the function $\phi:\mathbb{N}\rightarrow R^*$, defined by $\phi(i)=x^i$ is injective, which implies immediately that $R^*$ is an infinite set. The inverse of this function is the one defined in claim 1.

\item[5.]\label{Cc} From claim 4 it follows that any invertible element of a ring $R$ such that $R^*$ is finite has finite multiplicative order.  In addition, if $R^*=\{x_1,x_2,\dots,x_n\}$ and $M$ denotes the least common multiple of the set of natural numbers $\{o(x_1), o(x_2),\dots, o(x_n)\}$, then, 
\begin{equation}\label{GO}
x_i^M=1,
\end{equation}
for all $i\in \{1,2,\dots,n\}$.

\item[6.]\label{Dd}  The set of zero divisors of a ring does not have finite multiplicative order. In particular, nilpotent elements of a ring do not have finite multiplicative order.
\end{itemize}
A ring $R$ is said to have \textit{finite multiplicative order}, if there exists a natural number $M\geq 1$ such that $x^M=1$, for all $x\in R^*$. The smallest natural number $M$ such that $x^{M}=1$ for all $x\in R^*$ is called the \textit{multiplicative order} of the ring $R$ and it will be denoted by $o(R)$. If we set
\begin{equation*}\label{EulerPowers}
E(R)=\{\hspace{0.1cm}M\in\mathbb{N} \hspace{0.1cm}|\hspace{0.1cm} x^M=1 \text{ for all } x\in R^*\},
\end{equation*}
then to say that $R$ has finite multiplicative order is equivalent to saying that $E(R)$ is a nonempty set. If $R$ has finite multiplicative order then $o(R)=\min(E(R))$. 

\begin{itemize}

\item[7.]\label{Ee11} When $R$ has finite multiplicative order, from item 3 it follows that for all $M\in E(R)$ and for all $x\in R^*$, $o(x)|M$. Additionally, if it is assumed that there exists $y\in R^*$ with $o(y)=o(R)$, then it is clear that for all $M\in E(R)$, $o(R)|M$.  The condition $o(y)=o(R)$ holds, for instance when the group of units $R^*$ is a finite cyclic group.

\item[8.]\label{Ee} If $R^*=\{x_1,x_2,\dots,x_n\}$ is a finite group, claim 5 above implies that $R$ has finite multiplicative order. In addition, the multiplicative order of $R$ is the least common multiple of the set of natural numbers $\{o(x_1), o(x_2),\dots, o(x_n)\}$, i.e.,
\begin{equation}\label{EulerPM}
o(R)=\text{lcm}\{o(x_1),o(x_2),\dots,o(x_n)\}.
\end{equation}
In fact,  from item 3 above it follows that $M\in E(R)$, if and only if $M$ is a multiple of all the elements in the set $\{o(x_1),o(x_2),\dots,o(x_n)\}$. The latter certainly implies that $\min(E(R))=\text{lcm}\{o(x_1),o(x_2),\dots,o(x_n)\}$. Relation (\ref{EulerPM}) provides a clear method to compute $o(R)$: first, determine the group of units of the ring $R$; secondly, compute the multiplicative order of every element of the group $R^*$, and finally, determine the least common multiple of those numbers. Certainly, this method is not the best way to determine $o(R)$.

\item[9.]\label{Ff} If $R^*=\{x_1,x_2,\dots,x_n\}$ is a finite group, Lagrange theorem implies that
\begin{equation}\label{EulerPE}
x^{|R^*|}=1,\hspace{0.2cm}\text{for all}\hspace{0.2cm}x\in R^*,
\end{equation}
where $|R^*|$ denotes the order of $R^*$. In terms of the notation introduced above, Lagrange theorem shows that $|R^*|\in E(R)$, which implies that, if $R^*$ is a finite group, $|R^*|$ is an upper bound for $o(R)=\text{lcm}\{o(x_1),o(x_2),\dots,o(x_n)\}$.

\item[10.]\label{Gg} If $R=\mathbb{Z}_n$, Fermat-Euler theorem says that
\begin{equation}\label{EulerFer}
x^{\varphi(n)}= 1,\hspace{0.2cm}\text{for all}\hspace{0.2cm}x\in \mathbb{Z}_n^*,
\end{equation}
where $\varphi:\mathbb{N}\rightarrow\mathbb{N}$ is the Euler's totient function. In other words, $\varphi(n)$, counts the positive integers up to $n$ that are co-prime to $n$. The latter is equivalent to saying that $\varphi(n)=|\mathbb{Z}_n^*|$, thus, we can say that Fermat-Euler theorem is a consequence of Lagrange theorem. In terms of the notation introduced above, Fermat-Euler theorem implies that $\varphi(n)\in E(\mathbb{Z}_n)$.

\item[11.]\label{Hh} If $R=\mathbb{Z}_{12}$, $\mathbb{Z}_{12}^*=\{ 1,5,7,11 \}$, $\varphi(12)=4$ and every element of $\mathbb{Z}_{12}^*$ has multiplicative order 2.
So $|\mathbb{Z}_{12}^*|=4$ and $o(\mathbb{Z}_{12})=2$. This simple example shows that in general $o(R)<|R^*|$. 
\end{itemize}

From the discussion above, it is natural to say that a Fermat-Euler theorem holds in a ring $R$, if $E(R)$ is not an empty set and explicit formulas exist to describe elements in the set $E(R)$. When $R^*$ is finite, Lagrange theorem guarantees that $|R^*|\in E(R)$ (item 9 above), so in order to check that Fermat-Euler theorem holds in $R$, it is enough to prove that there exists an explicit formula to compute $M\in E(R)$, with  $o(R)\leq M\leq |R^*|$.\\

In this manuscript, it is shown that Fermat-Euler theorem holds for a general class of rings (Theorems \ref{Fermat} and \ref{Euler}). In most of our results, it is not necessary to assume that $R^*$ be finite nor $R$ commutative.


\section{Main results}
In this section, a generalization of Fermat-Euler theorem is established for a class of rings containing a collection of ideals satisfying the CNC-condition (Definition \ref{PosLiftIde}). For that purpose we need the following,
%
\begin{propo} \label{potencia}
Let $R$ be a ring and $N$ a nilpotent ideal of index $t \geq 2$ in $R$. Then the following statements hold:
\begin{enumerate}
	\item\label{pot} For any prime number $p$ such that $p \geq t$ and for all $n\in N$, $r\in R$ exists such that $$(1+n)^p = 1 + pnr.$$ 
	\item\label{port} Let $w$ be a natural number and $\bar{f}\in R/N$ be such that $(\overline{f})^w=\bar{1}$. If there exists a natural number $s>1$ such that $sN=0$, and all the prime factors of the number $s$ are greater than or equal to the nilpotency index $t$ of the ideal $N$, then 
$$
g^{ws}=1.
$$
for all $g\in\overline{f}$. In addition, if $w=o(\bar{f})$, then for all $g\in\bar{f}$, $o(g)$ divides $ws$. 
\end{enumerate}
\end{propo}

\proof 
\begin{enumerate}
\item  Since $n^t = 0$,
$$(1 + n)^p = \sum_{j=0}^{p}\binom{p}{j}n^j = 1 + \sum_{j=1}^{t-1}\binom{p}{j}n^j.$$
Since $p$ is a prime number, $p$ divides $\binom{p}{j}$ for all $1\leq j\leq p-1$.  Also, since $t\leq p$,
$$(1+n)^p=1+pn\left(k_1+{k_2}n+\cdots + k_{t-1}n^{t-2}\right)$$
where $k_i=\tbinom{p}{i}/p$. Therefore,
$$(1+n)^{p}= 1+pnr,$$
with $r=k_1+k_2n+\cdots + k_{t-1}n^{t-2}\in R$.

\item Let $p_1,p_2,p_3,\ldots ,p_m$ be the prime numbers, not necessarily distinct, appearing in the prime factor decomposition of the integer $s$, with $p_i \geq t$, for $i=1,2,3,\cdots,m.$ Since $(\overline{f})^w=(\overline{g})^w=\bar{1}$, there exists $n \in N$ such that $g^w = 1 + n$. Since $p_1\geq t$, from item \ref{pot}, there exists $r_1\in R$ such that
$$g^{wp_1}= (1 + n)^{p_1}= 1 + p_1nr_1.$$
Similarly, since $p_2\geq t$ and $p_1nr_1\in N$, item \ref{pot} implies that there exists $r_2\in R$ such that 
$$g^{wp_1p_2}=(1+p_1nr_1)^{p_2}=1+p_2(p_1nr_1)r_2.$$
Continuing this way, $r_3, r_4,\ldots, r_m\in R$ exists such that  
$$g^{ws}=1+sn(r_1r_2\cdots r_m).$$ 
In other words, 
$$g^{ws}=1+sh,$$ 
where $h=nr_1r_2\cdots r_m\in N$. Finally, since $h\in N$ and $sN=0$, it follows that $g^{ws}=1$. The proof of the last statement in claim \textit{2} is a consequence of claim 3, section \ref{BF}.
\end{enumerate}
\cqd

Note that under the assumptions of Proposition (\ref{potencia}),
the natural number $ws$ is not in general the multiplicative order of all the elements in the residue class $\overline{f}$. In fact, considering $R=\mathbb{Z}_{25}$ and $N=\langle 5 \rangle$, it is clear that $N$ has nilpotency  index $t=2$,  $sN={0}$ for $s=5$ and, the element $\overline{f}=\overline{4}=\{4,9,14,19,24\}$ in the quotient ring
$$
\frac{\mathbb{Z}_{25}}{\langle 5 \rangle}\cong \mathbb{Z}_5
$$
has multiplicative order $w=2$, that is $\overline{f}^2=\overline{1}$. So, from claim \ref{port} of Proposition (\ref{potencia}), it follows that any element $g\in \overline{f}$ is such that $g^{10}=1$. However, 10 is not the multiplicative order of all elements in $\overline{f}$. In fact, $24$ has multiplicative order $2$ in $\mathbb{Z}_{25}$, while the rest of the elements in $\overline{4}$ have multiplicative order equal to 10.\\

Note that the hypothesis in claim \ref{port} of Proposition \ref{potencia}, which requires that all prime factors of $s$ be greater or equal to the nilpotency index $t$ of the ideal $N$, restricts enormously the number of applications of that proposition. For instance, if we consider $R=\mathbb{Z}_{2^t}$ with $2\leq t$ and $N=\langle 2 \rangle$, it is clear that $N^{t}=\{0\}$  and $sN=\{0\}$ for $s=2^{t-1}$. Thus, according to claim \ref{port} of Proposition \ref{potencia}, in order to analyze the multiplicative order of an element $f$ in a ring $R$ by projecting it to the quotient ring $R/N\cong \mathbb{Z}_2$, it is necessary that $t\leq 2$, therefore $t=2$. Hence, the order of the elements in the ring $\mathbb{Z}_4$ can be described, which is easily done by hand. In the following lines, we show how to overcome such restrictions.\\


\begin{defin}{\cite[Definition 3.2]{liftidemp}}\label{PosLiftIde}
A collection $\{N_1,..., N_k\}$ of ideals of a ring $R$ satisfies the {\it CNC-condition} if the following properties hold:
\begin{enumerate}
\item \label{chc} {\bf Chain condition:} $\{0\}=N_k\subset N_{k-1}\subset\cdots \subset N_{2}\subset N_{1}\subset R$.
\item \label{nic} {\bf Nilpotency condition:} for $i=1,2,3,\ldots,k-1$, there exists $t_i \geq 2$ such that $N_i^{t_i}\subset N_{i+1}$.
\item \label{cac} {\bf Characteristic condition:} for $i=1,2,3,\ldots,k-1$, there exists $s_i\geq 1$ such that  $s_iN_i\subset N_{i+1}$. In addition, the prime factors of  $s_i$ are greater than or equal to $t_i$.   
\end{enumerate}
The minimum number $t_i$ satisfying the nilpotency condition will be called the nilpotency index of the ideal $N_i$ in the ideal $N_{i+1}$. Similarly,  
the minimum number $s_i$ satisfying the characteristic condition will be called the characteristic of the ideal $N_i$ in the ideal $N_{i+1}$. 
\end{defin}

\noindent
The nilpotency condition and the characteristic condition of the previous definition can be formulated as follows:

\begin{itemize}
\item[a.]\label{pronic} The nilpotency condition is equivalent to the following: for $i=1,2,\ldots,k-1$, $N_{i}/N_{i+1}$ is a nilpotent ideal of index $t_i$ in the ring $R/N_{i+1}$, (for details see {\cite[Definition 3.2]{liftidemp}}).  
\item[b.]\label{procac} The characteristic condition is equivalent to the following: for $i=1,2,\ldots,k-1$, there exists a natural number $s_i\geq1$ such that $s_i(N_{i}/N_{i+1})=0$ in the ring $R/N_{i+1}$, (for details see {\cite[Definition 3.2]{liftidemp}}).
\end{itemize}

The following theorem is one of the main results of this manuscript, allowing us to extend the classical  Fermat-Euler theorem to rings where the CNC-condition is satisfied.

\begin{theo}\label{IdemGeral} 
Let $	R$ be a ring, $\{N_1, N_2, \ldots, N_k\}$ a collection of ideals of $R$ satisfying the CNC-condition and $s_i$ the characteristic of the ideal $N_i$ in the ideal $N_{i+1}$.  Let $w$ be a natural number and $f\in R$ such that $(f+N_1)^{w}=1+N_1$. Then,
\begin{equation}\label{sat}
x^{ws_1s_2\cdots s_{k-1}}=1
\end{equation}
for all $x \in f+N_1\subset R$. In addition, if $w$ is the multiplicative order of $f+N_1$ in $R/N_1$, then the multiplicative order of any $x\in f+N_1$ must divide $ws_1s_2\cdots s_{k-1}$.
\end{theo} 

\proof Since $(f+N_1)^{w}=1+N_1$, for any $x \in f + N_1$, $(x+N_1)^{w}=1+N_1$. The isomorphism 
\begin{equation}\label{iso1}
R/N_1\cong \frac{(R/N_2)}{(N_1/N_2)},
\end{equation}
implies that $(x+N_2+N_1/N_2)^{w}=1+N_2+N_1/N_2$ in the ring $(R/N_2)/(N_1/N_2)$. Since $N_{1}/N_{2}$ is a nilpotent ideal of index $t_1$ in the ring $R/N_{2}$ and $s_1$ satisfies the hypothesis of claim 2 of Proposition \ref{potencia}, it follows that  
$$
(x+N_2)^{ws_1}=1+N_2,
$$
in the ring $R/N_2$. From the isomorphism
\begin{equation}\label{iso2}
R/N_2\cong \frac{(R/N_3)}{(N_2/N_3)},
\end{equation}
it follows that $(x+N_3+N_2/N_3)^{ws_1}=1+N_3+N_2/N_3$ in the ring $(R/N_3)/(N_2/N_3)$. Since $N_{2}/N_{3}$ is a nilpotent ideal of index $t_2$ in the ring $R/N_{3}$ and $s_2$ satisfies the hypothesis of claim 2  
of Proposition \ref{potencia}, it follows that 
$$
(x+N_3)^{ws_1s_2}=1+N_3,
$$
in the ring $R/N_3$. In the same way, since 
\begin{equation}\label{isogen}
R/N_i\cong \frac{(R/N_{i+1})}{(N_i/N_{i+1})},
\end{equation}
it follows that
$$
(x+N_{i+1})^{ws_1s_2\cdots s_{i}}=1+N_{i+1}.
$$
Finally, in the last step of the chain of ideals, 
$$
x^{ws_1s_2\cdots s_{k-1}}=1.
$$ 
The proof of the last claim in the theorem  is a consequence of item 3, Section \ref{BF}.
\cqd

\begin{obs}
It follows from the proof of the Theorem (\ref{IdemGeral}) that, if $x\in R$ is such that $w$ is the multiplicative order of $x+N_1$ in $R/N_1$, then the multiplicative order of the element $x+N_{2}$ in the ring $R/N_{2}$ must divide $ws_1$. Similarly, the multiplicative order of the element $x+N_{2}$ in the ring $R/N_{2}$ must divide $ws_1s_2$,  and so on. At the end of this process, it is seen that the multiplicative order of $x$ in $R$ must divide the natural number $ws_1s_2\cdots s_{k-1}$. The following chain of ring homomorphisms, 
$$
R\xrightarrow{\phi_{k-1}}\frac{R}{N_{k-1}}\xrightarrow{\phi_{k-2}} \cdots \xrightarrow{\phi_{3}}\frac{R}{N_{3}}\xrightarrow{\phi_{2}}\frac{R}{N_{2}}\xrightarrow{\phi_{1}}\frac{R}{N_{1}},
$$
appears naturally in that process. 
\end{obs}

We end this section with two important results that are a consequence of Theorem \ref{IdemGeral}.

\begin{theo}\label{Fermat}[An extension of Fermat's little Theorem]  
Let $	R$ be a ring, $\{N_1, N_2, \ldots, N_k\}$ a collection of ideals of $R$ satisfying the CNC-condition and $s_i$ the characteristic of the ideal $N_i$ in the ideal $N_{i+1}$. Then, the following statements hold:
\begin{enumerate}
\item\label{E1} If $(x+N_1)^{w}=1+N_1$, for all \hspace{0.1cm}$x+N_1\in (R/N_1)^*$, then for all $x\in R^*$, $x^{ws_1s_2\cdots s_{k-1}}=1$.
\item\label{E2} If $(R/N_1)^*$ is a finite group, then for all  $x\in R^*$, $x^{|(R/N_1)^*|s_1s_2\cdots s_{k-1}}=1$.
\item\label{E3} If $R^*$ is a finite group, then for all $x\in R^*$, $x^{|(R/N_1)^*||N_1|}=1$.
\end{enumerate}
\end{theo}

\proof For $x\in R^*$, it is clear that $x+N_1\in (R/N_1)^*$. By hypothesis $(x+N_1)^{w}=1+N_1$, so it follows from Theorem \ref{IdemGeral} that
\begin{equation}\label{M_1}
x^{ws_1s_2\cdots s_{k-1}}=1.
\end{equation}
To prove claim \textit{\ref{E2}}, as before, if $x\in R^*$ then $x+N_1\in (R/N_1)^*$. Since $(R/N_1)^*$ is a finite group, by Lagrange theorem $(x+N_1)^{|(R/N_1)^*|}=1+N_1.$
So from Theorem (\ref{IdemGeral}), it follows that 
\begin{equation}\label{M_2}
x^{|(R/N_1)^*|s_1s_2\cdots s_{k-1}}=1.
\end{equation}
For claim \textit{\ref{E3}}, since $R^*$ is a finite group, by Lagrange theorem $x^{|R^*|}=1$, for all $x\in R^*$. Since $R$ satisfies the CNC-condition, it is easy to see that $|R^*|=|(R/N_1)^*||N_1|$ (see alsoTheorem 3.3 in \cite{liftunits}). Therefore
\begin{equation}\label{M_3}
x^{|(R/N_1)^*||N_1|}=1,
\end{equation}
for all $x\in R^*$.
\cqd

In general terms, what has been proved in Theorem \ref{Fermat} is that if a Fermat-Euler theorem holds in $R/N_1$, then a Fermat-Euler theorem holds in $R$. Equivalently, under the assumptions of Theorem \ref{Fermat},
$$
M_1,M_2,M_3\in E(R),
$$
where,
\begin{equation}\label{DefMs}
M_1=ws_1s_2\cdots s_{k-1},\hspace{0.4cm}M_2=|(R/N_1)^*|s_1s_2\cdots s_{k-1},\hspace{0.2cm}\text{and}\hspace{0.2cm}M_3=|(R/N_1)^*||N_1|.
\end{equation}
If the quantities $M_1, M_2, M_3$ are assumed to be well defined (this is true when $|(R/N_1)^*|<\infty$ and $|N_1|<\infty$) and $w\leq |(R/N_1)^*|$ (this holds when $w=o(R/N_1)$), then it can be seen that 
\begin{equation}\label{des}
M_1\leq M_2\leq M_3.
\end{equation}
In fact, since $|N_1|=|N_1/N_2||N_2/N_3|\cdots|N_{k-1}/N_k|$ and $s_i(N_i/N_{i+1})=0\in R/N_{i+1}$ for all $i=1,2,\dots,k-1$, it follows that $s_i\leq |N_i/N_{i+1}|$, thus $M_2\leq M_3$. Under additional assumptions, the inequalities in (\ref{des}) can be improved, and
it is possible to show that $M_1|M_2$ and that $M_2|M_3$.\\ 

Next, an example that justifies the name of Theorem \ref{Fermat} is provided.

\begin{example}\label{Z_pi}\rm{
Let $p$ be a prime number, $k\in\mathbb{N}$ and let $R=\mathbb{Z}_{p^k}$ be the ring of integers modulo $p^{k}$. It is readily seen that $\{\langle p\rangle,\langle p^2\rangle,\dots,\langle p^k\rangle \}$ is a collection of ideals of $R$ satisfying the CNC-condition, with constant nilpotency indexes $t_i=2$ and constant characteristics $s_i=p$, for all $i=1,2,\dots, k-1$. Since
$\mathbb{Z}_{p^k}/\langle p\rangle\cong \mathbb{Z}_p$, it follows that $|\langle p\rangle|=p^{k-1}$, $|(\mathbb{Z}_{p^i}/\langle p\rangle)^*|=p-1$ and $o(\mathbb{Z}_{p^k}/\langle p\rangle)=p-1$. Therefore,  from items \textit{\ref{E1}}, \textit{\ref{E2}}, and  \textit{\ref{E3}} in Theorem (\ref{Fermat}), it follows that
$$
x^{(p-1)p^{i-1}}=1,
$$
for all $x\in (\mathbb{Z}_{p^k})^*=\{i\in\mathbb{Z}_{p^k}:i \not\equiv 0 \mod{(p)})\}$. In this example, the parameters $M_1, M_2, M_3$ defined in (\ref{DefMs}) are such that,
\begin{equation*}\label{igual}
M_1=M_2=M_3=\varphi(p^i).
\end{equation*}
}
\end{example}

In what follows several examples of rings where the parameters $M_1, M_2, M_3$ satisfy other relations will be provided. For example:
\begin{itemize}
\item $M_1=M_2<M_3$ (Example \ref{Horacio}).
\item $M_1<M_2<M_3$ (Example \ref{destricta}).
\item $M_3$ is not well defined ($|N_1|=\infty$) and $M_1= M_2$ (Example \ref{infinito}).
\end{itemize}

\begin{theo}\label{Euler}[A generalization of Euler's Theorem]  
Let $j\in\mathbb{N}$ and let $\{R_1,R_2,\dots,R_j\}$ be a collection of rings. Let $\{N_{i1}, N_{i2}, \ldots, N_{ik_i}\}$ be a collection of ideals of the ring $R_i$ satisfying the CNC-condition and $s_{il}$ the characteristic of the ideal $N_{il}$ in the ideal $N_{i(l+1)}$, $1\leq l\leq k_i-1$. If $(R_i/N_{i1})^*$ is a finite group  for all $1\leq i\leq j$, then 
$$
y^M=(1_{R_1},1_{R_2},\dots,1_{R_j}),
$$
for all $y=(y_1,y_2,\dots,y_j)\in (R_1)^*\times (R_2)^*\times\cdots\times (R_j)^*$, where $M$ is the smallest common multiple of all the numbers $\{m_1,m_2,\dots,m_j\}$ with $m_i=|(R_i/N_{i1})^*|s_{i1}s_{i2}\cdots s_{i(k_i-1)}$.
\end{theo}

\proof
Let $y=(y_1,y_2,\dots,y_j)\in (R_1)^*\times (R_2)^*\times\cdots\times (R_j)^*$. Since $M$ is the smallest common multiple of all the numbers in the set $\{m_1,m_2,\dots,m_j\}$,$r_1,r_2,\dots,r_j\in\mathbb{N}$ exist such that $M=m_jr_j$. Then
\begin{equation*}
y^M=\left(y_1^M,y_2^M,\dots,y_j^M\right)=\left((y_1^{m_1})^{r_1},(y_2^{m_2})^{r_2},\dots,(y_j^{m_j})^{r_j} \right).
\end{equation*}
From item \textit{\ref{E2}} of Theorem \ref{Fermat},  $y_i^{m_i}=1_{R_i}$, for all $1\leq i\leq j$. Therefore $y^M=(1_{R_1},1_{R_2},\dots,1_{R_j})$.
\cqd

It is worth mentioning that if the rings $R_i/N_{i1}$ are assumed have finite multiplicative order and the numbers $m_i=|(R_i/N_{i1})^*|s_{i1}s_{i2}\cdots s_{i(k_i-1)}$ appearing in the previous theorem are replaced by the numbers 
$$
m_i=w_is_{i1}s_{i2}\cdots s_{i(k_i-1)},\hspace{0.3cm}\text{with}\hspace{0.3cm} w_i=o(R_i/N_{i1}),
$$ 
a similar version of Theorem \ref{Euler} can be proved. Next, an example that justifies the name of Theorem \ref{Euler} is provided.

\begin{example} \rm{
Let $n$ be a natural number, and let $\mathbb{Z}_n$ be the ring of integers modulo $n$. The Chinese Remainder Theorem implies that
$$
(\mathbb{Z}_{n})^*\cong R_1^*\times R_2^*\times\cdots\times R_j^*
$$
where $R_i=\mathbb{Z}_{p_i^{k_i}}$ and $n=p_1^{k_1}p_2^{k_2}\cdots p_j^{k_j}$ is the prime decomposition of $n$. From example \ref{Z_pi}, it follows that 
$$
\{\langle p_i\rangle,\langle p_i^2\rangle,\dots,\langle p_i^{k_i}\rangle \}
$$
is a collection of ideals satisfying the CNC-condition in the ring $R_i$ with characteristics $s_{il}=p_i$, for $1\leq l\leq k_i-1$, and also that $|(R_i/\langle p_i\rangle)^*|=p_i-1$, for $1\leq i\leq j$. Thus, in terms of the notation of Theorem \ref{Euler}, $ m_i=(p_i-1)p_i^{k_i-1}$, $M=m_1m_2\cdots m_j$ and finally
$$
z^{(p_1-1)p_1^{k_1-1}(p_2-1)p_2^{k_2-1}\cdots(p_j-1)p_j^{k_j-1}}=1
$$
for all $z\in (\mathbb{Z}_n)^*=\{x\in\mathbb{Z}_n: x\not\equiv 0 \mod{(p_i)}, \text{ for all } 1\leq i\leq j \}$. Therefore, the classical Fermat-Euler theorem is recovered:
$$
z^{\varphi(n)}=1,\hspace{0.2cm}\text{for all}\hspace{0.2cm}z\in \mathbb{Z}_n^*.
$$
}
\end{example}


\section{Applications of the main results}

\quad In this section, Theorems \ref{Fermat} and \ref{Euler} will be applied in order to determine elements in the set $E(R)$ for several rings which include: rings containing a nilpotent ideal; matrix ring $M_n(R)$ where $R$ is a commutative ring containing a collection of ideals satisfying the CNC-condition; group ring $RG$ where the ring $R$ contains a nilpotent ideal; group rings $RG$ where $R$ is a chain ring; polynomial ring $R[x]$ where $R$ is a commutative ring containing a collection of ideals satisfying the CNC-condition. Examples are given illustrating the results.

\subsection{Rings containing a nilpotent ideal}

In the following result, for a ring $R$ containing a nilpotent ideal $N$, elements of $E(R)$ in terms of elements of the set $E(R/N)$ are determined.

\begin{propo}\label{GeralNil} 
Let $R$ be a ring and $N$ a nilpotent ideal of nilpotency index $k\geq 2$ in $R$.  Let $s>1$ be the characteristic of the quotient ring $R/N$. Then, the following statements hold:
\begin{enumerate}
\item\label{E1GN} If $(\overline{x})^{w}=\overline{1}$, for all \hspace{0.1cm}$\overline{x}\in (R/N)^*$, then for all $x\in R^*$, $x^{ws^{k-1}}=1$.
\item\label{E2GN} If $(R/N)^*$ is a finite group, then for all  $x\in R^*$, $x^{|(R/N)^*|s^{k-1}}=1$.
\item\label{E3GN} If $R^*$ is a finite group, then for all $x\in R^*$, $x^{|(R/N)^*||N|}=1$.
\end{enumerate}
\end{propo} 
\proof It is easy to show (see also \cite{liftunits} (Proposition 4.1) that the collection
$$
B=\{N, N^2,...,N^k\}
$$ 
of ideals of the ring $R$ satisfies the CNC-condition with nilpotency index and characteristic of the ideal $N^{i}$ in the ideal $N^{i+1}$ being $t_i=2$ and $s_i=s$ for all $i=1,2,3,\dots,k-1$. Therefore, the proof of this proposition is a clear application of Theorem \ref{Fermat}. 
\cqd

\begin{example}\label{Horacio}\rm{
Let $S$ be a commutative ring, $I$ a maximal nilpotent ideal of $S$ with nilpotency index $i\in\mathbb{N}$ and let $R=\{a+bu : a,b \in S, u^2=0 \}$. It is readily seen that $R$ with the (obvious) addition and multiplication operations is a commutative ring with cardinality $|R|=|S|^{2}$. It is also easily seen that $R$ is isomorphic to the ring of polynomials with coefficients in  $S$ modulo the ideal generated by $x^2$, that is $S[x]/\langle x^2\rangle$. It is easily checked that
$$
R^*=\{a+bu : a\in S^*, b\in S\},
$$
so the cardinality of $R^*$ is 
$$
|R^*|=|S^*||S|.
$$ Now, if $p$ denotes the characteristic of the residue field $\mathbb{F}=S/I$ and $p^t$ its cardinality, then, on the one hand, the Lagrange theorem implies that 
$|S|=p^t|I|$;  on the other hand, the proposition 4.1 in \cite{liftunits} implies that $|S^*|=(p^t-1)|I|.$
Thus, it is concluded that
\begin{equation}\label{Ho1}
x^{|R^*|}=
x^{p^t(p^t-1)|I|^2}=1,\hspace{0.5cm}\text{for all}\hspace{0.5cm}x\in R^*.
\end{equation}
Now, it is not difficult to check that $N_1=\{ a+bu : a\in I, b\in S\}$ is a maximal nilpotent ideal of the ring $R$ with nilpotency index $k=i+1$. In addition, 
$$
\frac{R}{N_1}\cong \frac{S}{I} =\mathbb{F},
$$  
hence the characteristic of the quotient ring $R/N_1$ is  $s=p$ and $w=o(R/N_1)=|(R/N_1)^*|=p^t-1$. Therefore, it is concluded from claims \ref{E1GNM} or \ref{E2GNM} of Proposition 4.1 that 
\begin{equation}\label{Ho2}
x^{(p^t-1)p^{i}}=1,\hspace{0.5cm}\text{for all}\hspace{0.5cm}x\in R^*.
\end{equation}
In this example, 
\begin{equation}\label{otraw}
M_1=M_2=(p^t-1)p^{i}, \hspace{0.3cm}\text{and}\hspace{0.3cm}M_3=p^{t}(p^t-1)|I|^2.
\end{equation}
When, for instance,  $S=\mathbb{Z}_{p^i}$ and $I=\langle p \rangle$ is the ideal generated by $p$ in $\mathbb{Z}_{p^i}$, it is not difficult to verify that $\mathbb{F}=\mathbb{Z}_p$, so it follows that $t=1$ and $|I|=p^{i-1}$. Therefore, it is deduced from (\ref{otraw}) that
$$
M_1=M_2=(p-1)p^{i}<M_3=(p-1)p^{2i-1}.
$$
 }
\end{example}

A consequence of the previous result is the following:

\begin{coro}\label{coroe}
Let $R$ be a commutative ring, $a \in R$ be a nilpotent element of index $k$ and $s>1$ be the characteristic of the quotient ring $R/\langle a\rangle$. Then, the following statements hold
\begin{enumerate}
\item\label{E1EN} If $(\overline{x})^{w}=\overline{1}$, for all \hspace{0.1cm}$\overline{x}\in (R/\langle a\rangle)^*$, then for all $x\in R^*$, $x^{ws^{k-1}}=1$.
\item\label{E2EN} If $(R/\langle a\rangle)^*$ is a finite group, then for all  $x\in R^*$, $x^{|(R/\langle a\rangle)^*|s^{k-1}}=1$.
\item\label{E3EN} If $R^*$ is a finite group, then for all $x\in R^*$, $x^{|(R/\langle a\rangle)^*||\langle a\rangle|}=1$.
\end{enumerate}
\end{coro} 
\proof Since $R$ is a commutative ring, $\langle a \rangle$ is a nilpotent ideal of nilpotency index $k$ in $R$, and the result follows immediately from Proposition \ref{GeralNil}
\medskip
\cqd


\subsection{Matrix ring}

The following results determine elements of the set $E(M_n(R))$ in terms of the  elements of $E(M_n(R/N_1))$, where $R$ is a commutative ring containing a collection of ideals $\{N_1,\dots,N_k\}$ satisfying the CNC-condition.  
\begin{propo}\label{MatrixRing} 
Let $R$ be a commutative ring and let $M_n(R)$ be the ring of  $n\times n$ matrices with entries in $R$. Let $\{N_1, N_2, \ldots, N_k\}$ be a collection of ideals of $R$ satisfying the CNC-condition, and  
let  $s_i$ be the characteristic of the ideal $N_i$ in the ideal $N_{i+1}$. Then, the following statements hold:
\begin{enumerate}
\item\label{E1MR} If $(\overline{x})^{w}=\overline{1}$, for all \hspace{0.1cm}$\overline{x}\in (M_n(R/N_1))^{*}$, then for all $x\in (M_n(R))^*$, $x^{ws_1s_2\cdots s_{k-1}}=1$.
\item\label{E2MR} If $(M_n(R/N_1))^{*}$ is a finite group, then for all  $x\in (M_n(R))^*$, $x^{|(M_n(R/N_1))^{*}|s_1s_2\cdots s_{k-1}}=1$.
\item\label{E3MR} If $(M_n(R))^{\ast}$ is a finite group, then for all $x\in (M_n(R))^*$, $x^{|(M_n(R/N_1))^{\ast}||N_1|^{n^2}}=1$.
\end{enumerate}
\end{propo} 
\proof It is easy to show (see also \cite{liftunits} (Proposition 4.5)) that the collection  
$$
B=\{M_n(N_1), M_n(N_2),...,M_n(N_k)\}
$$ 
of ideals of the ring $M_n(R)$ satisfies the CNC-condition, with nilpotency index and characteristic of the ideal $M_n(N_i)$ in the ideal $M_n(N_{i+1})$ being exactly the same nilpotency index and characteristic of the ideal $N_i$ in the ideal $N_{i+1}$. In addition; it is also easy to prove (see also \cite{liftunits} (Proposition 4.5)) that,
$$
|(M_n(R))^{\ast}|=|(M_n(R/N_1))^{\ast}||N_1|^{n^2}.
$$
Therefore, the proof of this proposition follows from Theorem \ref{Fermat}.
\cqd

\begin{coro}\label{GeralNilMnR} 
Let $R$ be a commutative ring and let $M_n(R)$ be the ring of  $n\times n$ matrices with entries from $R$. Let $N$ be a nilpotent ideal of $R$ of index $k$ in $R$, and  let  $s$ be the characteristic of the quotient ring $R/N$.
Then, the following statements hold:
\begin{enumerate}
\item\label{E1GNM} If $(\overline{x})^{w}=\overline{1}$, for all \hspace{0.1cm}$\overline{x}\in (M_n(R/N))^{*}$, then for all $x\in (M_n(R))^*$, $x^{ws^{k-1}}=1$.
\item\label{E2GNM} If $(M_n(R/N))^{*}$ is a finite group, then for all  $x\in (M_n(R))^*$, $x^{|(M_n(R/N))^{*}|s^{k-1}}=1$.
\item\label{E3GNM} If $(M_n(R))^{\ast}$ is a finite group, then for all $x\in (M_n(R))^*$, $x^{|(M_n(R/N))^{\ast}||N|^{n^2}}=1$.
\end{enumerate}
\end{coro}
\proof 
As mentioned in the proof of proposition \ref{GeralNil}, the collection $\{N, N^2,...,N^{k}\}$ of ideals of the ring $R$ satisfies the CNC-condition with constant characteristic $s_i=s$, for all $i=1,2,3,\cdots, k-1$. Thus, the proof of this result is an immediate consequence of Proposition \ref{MatrixRing}.
\cqd 

\begin{coro}\label{GeralEleNilMnR} 
Let $R$ be a commutative ring and $M_n(R)$ be the ring of  $n\times n$ matrices with entries from $R$. Let $a$ be a nilpotent element  of  index $k$ in $R$ and  $s$ be the characteristic of the quotient ring $R/\langle a\rangle$. Then, the following statements hold:
\begin{enumerate}
\item\label{E1GNM} If $(\overline{x})^{w}=\overline{1}$, for all \hspace{0.1cm}$\overline{x}\in (M_n(R/\langle a\rangle))^{*}$, then for all $x\in (M_n(R))^*$, $x^{ws^{k-1}}=1$.
\item\label{E2GNM} If $(M_n(R/\langle a\rangle))^{*}$ is a finite group, then for all  $x\in (M_n(R))^*$, $x^{|(M_n(R/\langle a\rangle))^{*}|s^{k-1}}=1$.
\item\label{E3GNM} If $(M_n(R))^{\ast}$ is a finite group, then for all $x\in (M_n(R))^*$, $x^{|(M_n(R/\langle a\rangle))^{\ast}||\langle a\rangle|^{n^2}}=1$.
\end{enumerate}

\end{coro}
\proof It is enough to observe that the ideal $N=\langle a\rangle$ is nilpotent of index $k$ in the ring $R$. The result follows from Corollary \ref{GeralNilMnR}.
\cqd 
%

\begin{example}\rm{
Let $R$ be a commutative ring, $N$ a maximal nilpotent ideal of $R$ with nilpotency index $i\in\mathbb{N}$ and such that the residue field $R/N=\mathbb{F}$ has $p^t$ elements, where $p$  is a prime number denoting the characteristic of the field $\mathbb{F}$ and $t$ is a natural number.
Let $M_2(R)$ be the ring of $2\times 2$ matrices with entries from the ring $R$. Since, by assumption $N$ is a nilpotent ideal of index $i$ in $R$, corollary \ref{GeralNilMnR} can be applied. Since 
$$
|(M_2(\mathbb{F}))^*|=(|\mathbb{F}|^2-1)(|\mathbb{F}|^2-|\mathbb{F}|)=(p^{2t}-1)(p^{2t}-p^{t})
$$
and the characteristic of $\mathbb{F}$ is $s=p$, the item \ref{E2GNM} of the corollary implies that
\begin{equation}\label{Ma1}
x^{(p^{2t}-1)(p^{2t}-p^t)p^{i-1}}=1,\hspace{0.5cm}\text{for all}\hspace{0.5cm}x\in (M_2(R))^*.
\end{equation}
On the other hand, it follows from item \ref{E3GNM} of the same corollary that
\begin{equation}\label{Ma2}
x^{(p^{2t}-1)(p^{2t}-p^t)|N|^{4}}=1,\hspace{0.5cm}\text{for all}\hspace{0.5cm}x\in (M_2(R))^*.
\end{equation}
Observe that in this example 
\begin{equation}\label{oty}
M_2=(p^{2t}-1)(p^{2t}-p^t)p^{i-1}\hspace{0.4cm}\text{and}\hspace{0.4cm} M_3=(p^{2t}-1)(p^{2t}-p^t)|N|^{4}.
\end{equation}
When, for instance,  $R=\mathbb{Z}_{p^i}$ and $N=\langle p \rangle$ is the ideal generated by $p$ in $\mathbb{Z}_{p^i}$, it is not difficult to verify that $\mathbb{F}=\mathbb{Z}_p$, so it follows that $t=1$ and $|N|=p^{i-1}$. Therefore, it is deduced from (\ref{oty}) that
$$
M_2=(p^{2}-1)(p^{2}-p)p^{i-1}<M_3=(p^{2}-1)(p^{2}-p)p^{4(i-1)}.
$$
 }
 \end{example}


\subsection{Group rings}

If $R$ is a ring containing a collection of ideals satisfying the CNC-condition and $G$ is a group, by invoking Theorem \ref{Fermat}, some elements of the set $E(RG)$ are determined. 

\begin{propo}\label{GeralGR} 
Let $	R$ be a ring and $G$ a group. Let $\{N_1, N_2, \ldots, N_k\}$ be a collection of ideals of $R$ satisfying the CNC-condition with $s_i$ being the characteristic of the ideal $N_i$ in the ideal $N_{i+1}$. Then, the following statements hold:
\begin{enumerate}
\item\label{E1GR} If $(\overline{x})^{w}=\overline{1}$, for all \hspace{0.1cm}$\overline{x}\in ((R/N_1)G)^{*}$, then for all $x\in (RG)^*$, $x^{ws_1s_2\cdots s_{k-1}}=1$.
\item\label{E2GR} If $((R/N_1)G)^{*}$ is a finite group, then for all  $x\in (RG)^*$, $x^{|((R/N_1)G)^{*}|s_1s_2\cdots s_{k-1}}=1$.
\item\label{E3GR} If $(RG)^{\ast}$ is a finite group, then for all $x\in (RG)^*$, $x^{|((R/N_1)G)^{\ast}||N_1|^{|G|}}=1$.
\end{enumerate}
\end{propo} 

\proof 
Again it is easy to prove (see also \cite{liftunits} (Proposition 4.9)) that the collection  
$$
B=\{N_1G, N_2G, \ldots, N_kG\}
$$ 
of ideals of the ring $RG$ satisfies the CNC-condition with nilpotency index and characteristic of the ideal $N_iG$ in the ideal $N_{i+1}G$ being exactly the same nilpotency index and characteristic of the ideal $N_i$ in the ideal $N_{i+1}$. Furthermore it also can be seen that
$$
|(RG)^*|=|((R/N_1)G)^{\ast}||N_1|^{|G|}.
$$
Therefore, the proof of this proposition follows at once from Theorem \ref{Fermat}.
\cqd

\begin{coro}\label{GeralNilGR} 
Let $G$ be a group,  let $R$ be a ring, $N$ a nilpotent ideal of index $k$ in $R$, and $s$ the characteristic of the quotient ring $R/N$. Then, the following statements hold:
\begin{enumerate}
\item\label{E1GRN} If $(\overline{x})^{w}=\overline{1}$, for all \hspace{0.1cm}$\overline{x}\in ((R/N)G)^{*}$, then for all $x\in (RG)^*$, $x^{ws^{k-1}}=1$.
\item\label{E2GRN} If $((R/N)G)^{*}$ is a finite group, then for all  $x\in (RG)^*$, $x^{|((R/N)G)^{*}|s^{k-1}}=1$.
\item\label{E3GRN} If $(RG)^{\ast}$ is a finite group, then for all $x\in (RG)^*$, $x^{|((R/N)G)^{\ast}||N|^{|G|}}=1$.
\end{enumerate}
\end{coro}
\proof 
The proof is a direct consequence of Proposition \ref{GeralGR} and the fact that the collection $\{N, N^2,...,N^{k}\}$ of ideals of the ring $R$ satisfies the CNC-condition with constant characteristic $s_i=s$ for all $i=1,2,3,\cdots, k-1$.
\cqd 

\begin{coro}\label{cyu}
Let $	R$ be a commutative ring, let $a$ be a nilpotent element of index $k$ in $R$, $G$ a group and $s$ the characteristic of the quotient ring $R/\langle a\rangle$. Then, the following statements hold:
\begin{enumerate}
\item\label{E1GREN} If $(\overline{x})^{w}=\overline{1}$, for all \hspace{0.1cm}$\overline{x}\in ((R/\langle a\rangle)G)^{*}$, then for all $x\in (RG)^*$, $x^{ws^{k-1}}=1$.
\item\label{E2GREN} If $((R/\langle a\rangle)G)^{*}$ is a finite group, then for all  $x\in (RG)^*$, $x^{|((R/\langle a\rangle)G)^{*}|s^{k-1}}=1$.
\item\label{E3GREN} If $(RG)^{\ast}$ is a finite group, then for all $x\in (RG)^*$, $x^{|((R/\langle a\rangle)G)^{\ast}||\langle a\rangle|^{|G|}}=1$.
\end{enumerate}
\end{coro} 
\proof 
Since $R$ is a commutative ring, the ideal $N = \langle a \rangle $ is nilpotent of index $k$ in $R$, and the result follows from Corollary \ref{GeralNilGR}.
\cqd
\begin{example}\label{destricta}\rm{
Let $p$ be a prime number, $i\in\mathbb{N}$ and let $R=\{a+bg : a,b\in \mathbb{Z}_{p^i}, g^2=1 \}$ be the group ring $\mathbb{Z}_{p^i}G$ where $G=\{1,g\}$ is the cyclic group of order $n=2$. It is readily seen that $R$ with the (obvious) addition and multiplication operations is a commutative ring with cardinality $|R|=p^{2i}$. It is also evident that $R$ is isomorphic to the ring of polynomials with coefficients in  $\mathbb{Z}_{p^i}$ modulo the ideal generated by $x^2-1$ in $R$, that is $\mathbb{Z}_{p^i}G\cong \mathbb{Z}_{p^i}[x]/\langle x^2-1\rangle$. An easy calculation shows that
$$
(\mathbb{Z}_pG)^*=\{a+bu : a\neq b, a\neq -b\},
$$
so the cardinality of $(\mathbb{Z}_pG)^*$ is $|(\mathbb{Z}_pG)^*|=(p-1)^2$, when $p>2$ and  $|(\mathbb{Z}_2G)^*|=2$. Also, since $\langle p \rangle$ has nilpotency index $k=i$ in $\mathbb{Z}_{p^i}$, and 
$$
\frac{Z_{p^i}G}{\langle p\rangle G}\cong \mathbb{Z}_pG,
$$
then $|(\mathbb{Z}_{p^i}G)^*|=(p-1)^2p^{2(i-1)}$, when $p>2$ and  $|(\mathbb{Z}_{2^i}G)^*|=2^{2i-1}$. Thus, from item \textit{\ref{E3}} of Theorem (\ref{Fermat}), it follows that
for all $x\in R^*$
\begin{equation}\label{Ce1}
\left\{
\begin{aligned}
x^{2^{2i-1}}&=1,\hspace{0.7cm}\text{ if } p=2,\\
x^{(p-1)^2p^{2(i-1)}}&=1,\hspace{0.7cm}\text{ if } p>2. 
\end{aligned}
\right.
\end{equation}
On the other hand, from item \textit{\ref{E2}} in Theorem (\ref{Fermat}),
it follows that 
\begin{equation}\label{Ce1}
\left\{
\begin{aligned}
x^{2^{i}}&=1,\hspace{0.7cm}\text{ if } p=2,\\
x^{(p-1)^2p^{i-1}}&=1,\hspace{0.7cm}\text{ if } p>2, 
\end{aligned}
\right.
\end{equation}
for all $x\in R^*$. Finally, for all $x\in \mathbb{Z}_pG$, $x^{p-1}=1$, from  item \textit{\ref{E1}} of Theorem (\ref{Fermat}) with $w=p-1$,
\begin{equation}\label{Ce2}
\left\{
\begin{aligned}
x^{2^{i}}&=1,\hspace{0.7cm}\text{ if } p=2,\\
x^{(p-1)p^{i-1}}&=1,\hspace{0.7cm}\text{ if } p>2. 
\end{aligned}
\right.
\end{equation}
for all $x\in R^*$. In this example, when $p>2$, $M_1<M_2<M_3$.
}
\end{example}


\subsection{Commutative group rings $RG$ with $R$ a chain ring}

Let $R$ be a commutative chain ring and $G$ a commutative group. It is well known that $R$ contains a unique maximal nilpotent ideal $N=\langle a\rangle$ for some $a\in R$. If $k$ denotes the nilpotency index of $a$, and $p$ denotes the characteristic of the residue field $\F=R/\langle a\rangle$, from Corollary \ref{cyu}, it follows that:  
\begin{itemize}
\item If $(\F G)^{*}$ is a finite group, then for all  $x\in (\F G)^*$, $x^{|(\F G)^{*}|p^{k-1}}=1$.
\item If $(RG)^{*}$ is a finite group, then for all $x\in (R G)^*$, $x^{|(\F G)^{\ast}||\langle a\rangle|^{|G|}}=1$.
\end{itemize}
Examples of finite commutative chain rings include the ring of modular integers $R=\mathbb{Z}_{p^k}$, where $p$ is a prime number and $k>1$ is an integer. In this example the maximal nilpotent ideal is $N=\langle p \rangle$, with nilpotency index equal to $k$ in $\mathbb{Z}_{p^k}$, and $\F\cong \Z_p$. Thus, if $(RG)^{*}$ is a finite group,
\begin{equation}\label{few1}
x^{|(\mathbb{Z}_pG)^{*}|p^{k-1}}=1,\hspace{0.5cm}\text{ for all }\hspace{0.5cm} x\in (\mathbb{Z}_{p^k}G)^{*}.
\end{equation}
 
Another interesting example of chain rings is represented by Galois rings, $R=\mathbb{Z}_{p^k}[x]/\langle q(x) \rangle$ where $p$ is a prime number, $k>1$ and $q(x)$ is a monic polynomial of degree $r$ whose image in $\Z_{p}[x]$ is irreducible. In this example the maximal nilpotent ideal is $N=pR$, with nilpotency index equal to $k$ in $R$, and $\F=R/N\cong \F_{p^{r}}$. Thus, 
\begin{equation}\label{few2}
g^{|(\F_{p^r}G)^{*}|p^{k-1}}=1,\hspace{0.5cm}\text{ for all }\hspace{0.5cm} g\in ((\mathbb{Z}_{p^k}[x]/\langle q(x) \rangle)G)^{\ast}.
\end{equation}
In particular, if $G={e}$, the trivial group, it follows that
\begin{equation}\label{few2}
g^{(p^r-1)p^{k-1}}=1,\hspace{0.5cm}\text{ for all }\hspace{0.5cm} g\in (\mathbb{Z}_{p^k}[x]/\langle q(x) \rangle)^{\ast}.
\end{equation}


\subsection{Polynomial rings}

It is not difficult to prove that if $R$ is an associative ring with identity such that $R$ contains a family of ideals $\{N_1,N_2,\dots,N_k\}$ satisfying the CNC-condition with characteristics $s_1,s_2,\dots,s_{k-1}$ and nilpotency indexes $t_1,t_2,\dots,t_{k-1}$, then the family 
$$
\{N_1[x],N_2[x],\dots,N_{k-1}[x] \}
$$ 	
of ideals satisfies the CNC-conditions in the polynomial ring $R[x]$. Thus,  again from Theorem \ref{Fermat}, it follows that if $(g+N_1[x])^w=1+N_1[x]$ for all $g+N_1[x]\in (R[x]/N_1[x])^*$, then 
\begin{equation}\label{poly}
g^{ws_1s_2\cdots s_{k-1}}=1, \hspace{0.5cm}\text{ for all }\hspace{0.5cm} g\in (R[x])^*.
\end{equation}

\begin{example}\label{infinito}
\rm{Let $p$ be a primer number, $k\in\mathbb{N}$ and $R=\mathbb Z_{p^{k}}[x]$ the ring of polynomials with coefficients in the ring of integers modulo $p^{k}$. It is clear that $N_1=\langle p \rangle$ has nilpotency index $k$ in $\mathbb{Z}_{p^k}$, so the family of ideals $\{\langle p\rangle,\langle p^2\rangle,\dots, \langle p^k\rangle\}$ satisfies the CNC-condition in $\mathbb{Z}_{p^k}$ with constant characteristics $s_i=p$ and constant nilpotency indexes $t_i=2$. In addition, 
$$
\frac{\mathbb{Z}_{p^k}[x]}{N_1[x]}\cong \mathbb{Z}_{p}[x],
$$
so $(\mathbb{Z}_{p^k}[x]/N_1[x])^*\cong(\mathbb{Z}_{p})^*$. The latter implies that $o(\mathbb{Z}_{p^k}[x]/N_1[x])=p-1=w$. Therefore, it is concluded from (\ref{poly}) that
$$
g^{(p-1)p^{i-1}}=1,\hspace{0.5cm}\text{for all}\hspace{0.5cm}g\in R^*. 
$$
Note that in this example, $M_3$ is not well defined ($|N_1[x]|=\infty$) and $M_1=M_2$.}
\end{example}


\end{document}